\documentclass[reqno]{amsart}
\usepackage{hyperref,graphicx}

\begin{document}
\title[Cone Calculus]
{Introduction to the Analysis on Manifolds with Conical Singularities}

\author[Elmar Schrohe\hfil \hfilneg] {Elmar Schrohe}  

\address{Elmar Schrohe, Leibniz University Hannover,
Institute of Analysis, \mbox{Welfengarten 1}, 30167 Hannover, Germany}
\email{schrohe@math.uni-hannover.de}


\begin{abstract}
In these notes I will recall central elements of the cone calculus. The focus lies on conically degenerate differential operators and the Laplace-Beltrami operator with respect to a conically degenerate metric as a prototypical example. 
We will get to know manifolds with conical singularities, the Mellin transform, cone Sobolev spaces, and the notion of ellipticity in terms of the invertibility of the principal pseudodifferential symbol and the principal Mellin symbol. Finally I will sketch the full cone calculus. 
\end{abstract}

\maketitle \numberwithin{equation}{section}
\newtheorem{theorem}{Theorem}[section]
\newtheorem{corollary}[theorem]{Corollary}
\newtheorem{lemma}[theorem]{Lemma}
\newtheorem{remark}[theorem]{Remark}
\newtheorem{problem}[theorem]{Problem}
\newtheorem{example}[theorem]{Example}
\newtheorem{definition}[theorem]{Definition}
\allowdisplaybreaks

{

\section{Introduction}
We start with an example that illuminates some of the features that come up in connection with conical singularities. We consider a closed angular domain $B$ in the plane ${\mathbb R}^2$ with coordinates $x$ and $y$. We assume that the origin is the vertex of the domain and one of the bounding half lines is the non-negative $x$-axis while the other intersects the $x$-axis under an angle $0<{\alpha}<2\pi$. 

We are interested in solutions to the homogeneous Dirichlet problem 
\begin{eqnarray}\Delta u=0&&\text{ in } B^\circ\label{e.1.1}\\
 u=0 &&\text{ on } \partial B\setminus \{0\}\label{e.1.2}.
\end{eqnarray}
Naively, one might expect the problem to only have the trivial solution $u\equiv0$.
We will see, however, that something unexpected can happen.

We introduce polar coordinates $(t,{\theta})$ at the origin, so that $B$ can be written  
$$B=\{(t\cos {\theta},t\sin{\theta})): t\ge 0, 0\le {\theta}\le {\alpha}\}.$$
We let $v(t,{\theta}) = u(t\cos {\theta}, t\sin{\theta})$. 
Equations \eqref{e.1.1} and \eqref{e.1.2} become
\begin{eqnarray}
((t\partial_t)^2+\partial^2_{\theta} )v(t,{\theta})=0&&\text{in } (0,\infty)\times(0,{\alpha}),\label{e.1.3}\\
v(t,{\theta})=0&&\text{for } t>0, {\theta}=0,{\alpha}\label{e.1.4}  .
\end{eqnarray}
By ${\mathcal M}$  we denote the Mellin transform, which we apply with respect to $t$:
$${\mathcal M} v(z,{\theta}) = \int_0^\infty t^z v(t,{\theta}) \frac{dt}t, \ z\in {\mathbb C}.$$
Ignoring for the moment problems of convergence of the integral we see  that $z{\mathcal M} v(z) = {\mathcal M}((-t\partial_t)v)(z).$
Then Equations \eqref{e.1.3} and  \eqref{e.1.4} become
\begin{eqnarray}
(z^2+\partial_{\theta}^2){\mathcal M} v(z,{\theta})=0,&& z\in {\mathbb C}, {\theta} \in (0,{\alpha}) \label{e.1.5}\\
{\mathcal M}  v(z,{\theta})=0&&\text{for } {\theta}=0,{\alpha}. \label{e.1.6}
\end{eqnarray}
Now we note that the ODE 
\begin{eqnarray*}
(z^2+\partial_{\theta}^2)w({\theta})=0,&& z\in {\mathbb C}, {\theta} \in  (0,{\alpha})\label{e.1.7}\\
w({\theta})=0&&\text{for } {\theta}=0,{\alpha}. \label{e.1.8}
\end{eqnarray*}
has nonzero solutions for $z=j\pi/\alpha$, namely $\sin(j\pi{\theta}/{\alpha})$ for 
$j\in\mathbb Z\setminus\{0\}$, 
and the functions
$$t^{j\pi/\alpha}\sin(j\pi{\theta}/{\alpha})
$$
actually solve the problem \eqref{e.1.1}, \eqref{e.1.2}.  
Note that $t^{j\pi/{\alpha}}$ gets more and more singular as $j\to -\infty$.
\begin{center}
\resizebox*{0,5\textwidth}{!}{\includegraphics*{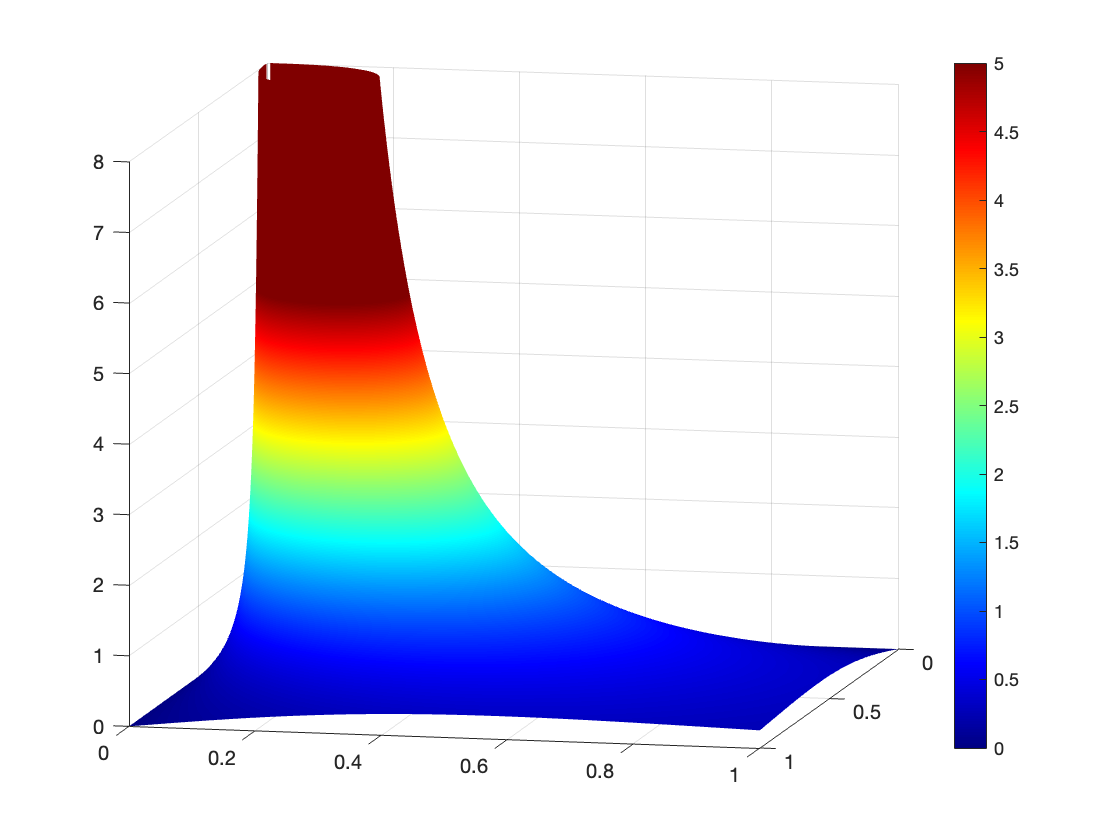}}\\
Figure 1: The singular solution $t^{-2} \sin(2\theta)$.
\end{center}
The first unbounded solution is $u_1(t,{\theta}) = t^{-\pi/{\alpha}}\sin(\pi{\theta}/{\alpha})$. It is in $L^2$ near $t=0$ if and only if ${\alpha}>\pi$, i.e. the angle is obtuse.
The presence of such singular solutions may create problems in numerical computations, so it is important to know them and to take them into account. 

\section{Manifolds with Conical Singularities}

\begin{definition}
A manifold with conical singularities $B$ is a topological space which is a smooth manifold outside a finite number of points while, close to each of these points, it has the structure of a cone whose cross-section is a smooth $n$-dimensional manifold $X$, i.e., there is a neighborhood  $U$ of this point in $B$ with 
$$U\cong X\times [0,1)/X\times \{0\}.$$
Moreover, $B$ comes with a smooth structure that preserves the cones. 
\end{definition}  	
For more details see \cite[Section 1]{SS1994} or  \cite[Section 3]{LS2019}. 

It makes the analysis much easier to blow up the tip and to work on the blown-up space which is a smooth manifold ${\mathbb B}$ of dimension $n+1$ with boundary $\partial {\mathbb B}=X$. We know that every manifold ${\mathbb B}$ with boundary has a `collar neighborhood' of the boundary, i.e. a neighborhood of $\partial {\mathbb B}$ of the form $[0,1)\times \partial {\mathbb B}$. In this neighborhood we choose local coordinates $(t,x)$ in $[0,1)\times X$, so that $t$ `measures' the distance to the boundary. 
\begin{center}
\resizebox*{0,3\textwidth}{!}{\includegraphics*{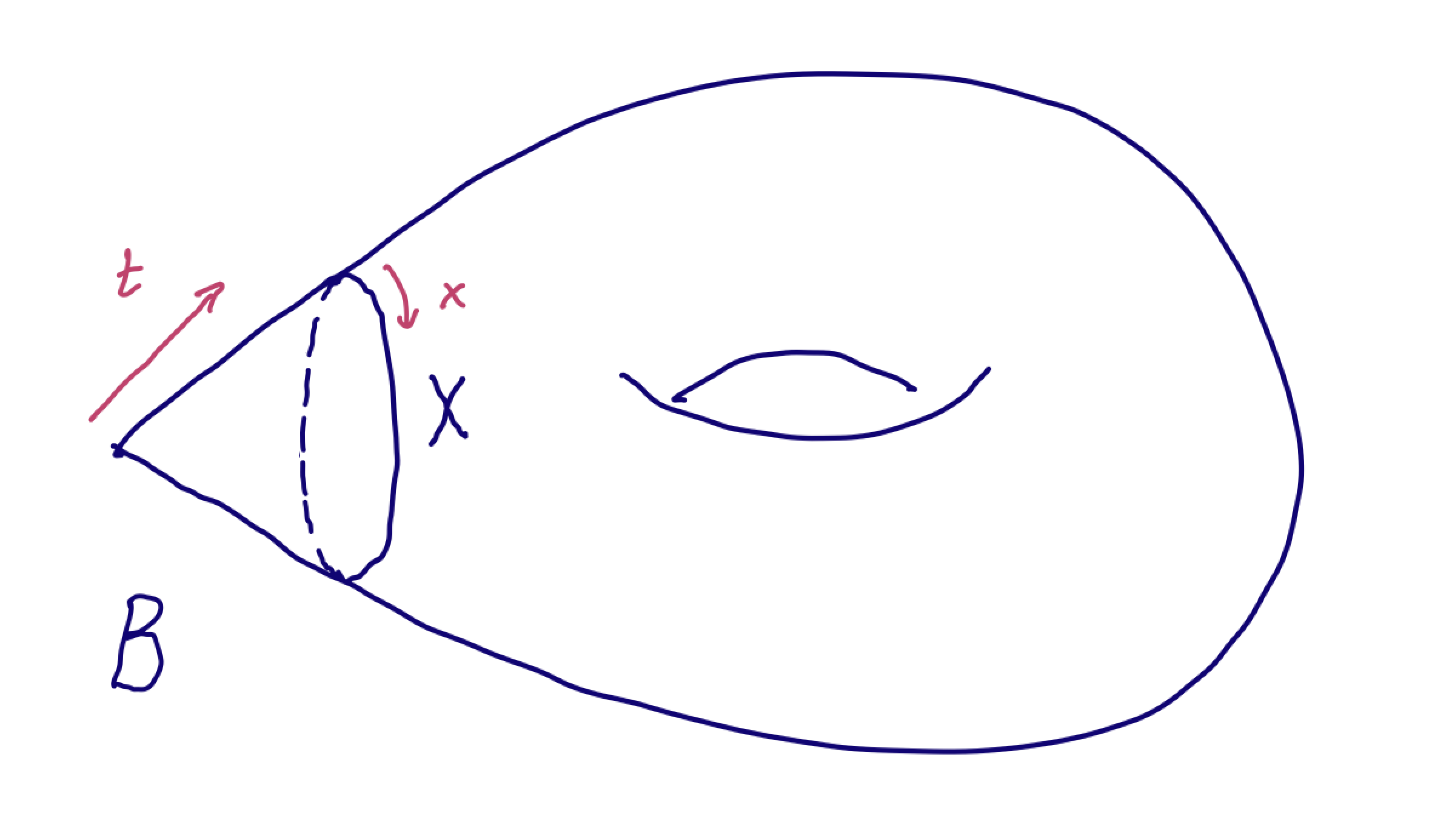}}\resizebox*{0,1\textwidth}{!}{\includegraphics*{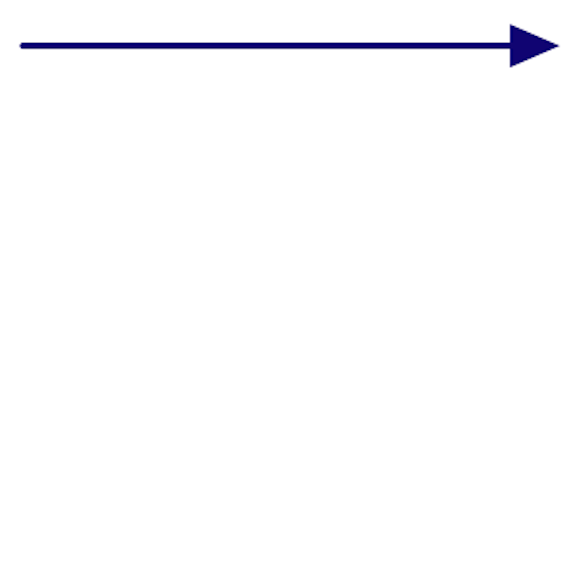}}\resizebox*{0,3\textwidth}{!}{\includegraphics*{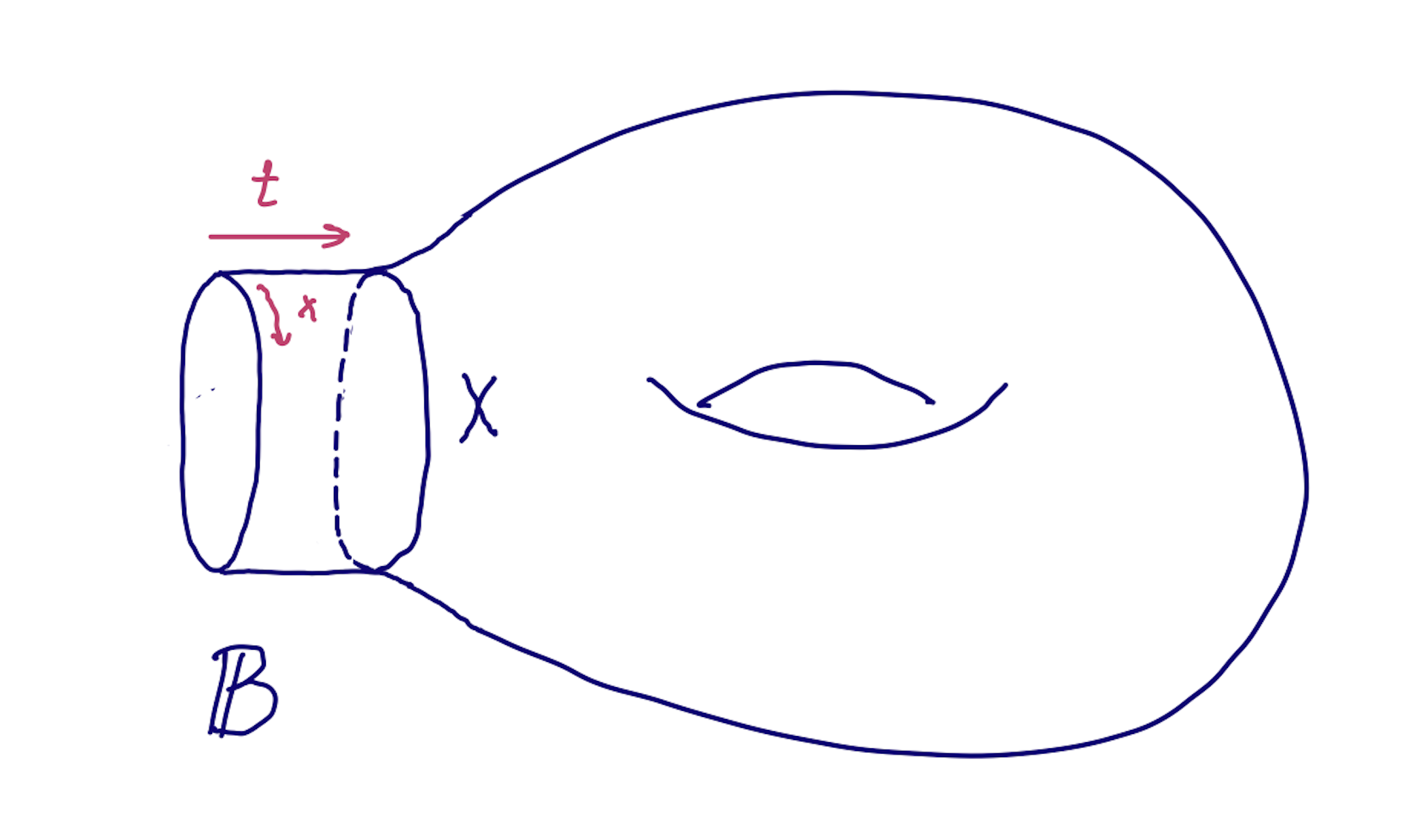}}\\
Figure 2. A manifold with a conical singularity $B$ and its blow-up ${\mathbb B}$
\end{center}  
We model the geometry of $B$ on ${\mathbb B}$ by endowing ${\mathbb B}$ with a  Riemannian metric $g$ that degenerates near the boundary, namely
\begin{eqnarray}
\label{eq.g}
g= dt^2 + t^2h(t), 
\end{eqnarray}
where $t\mapsto h(t)$ is a smooth family of non-degenerate Riemannian metrics on the cross-section $X$. One speaks of a straight cone, if $h$ is independent of $t$ and of a warped cone, if $h$ varies with $t$. Recall that \eqref{eq.g} is another way of representing the metric $g$ in local coordinates $(t, x)$  as the bilinear form given by the matrix 
\begin{eqnarray}
\label{eq.matrixg}
g= 
\begin{pmatrix}
1&0 \ldots 0 \\
0& \\
\vdots &t^2h(t)\\
0&   
\end{pmatrix} .
\end{eqnarray}
 
Notice that $g$ models precisely the conical situation:  Distances to the boundary remain unchanged (the $dt^2$ term) while - in the situation of a straight cone - distances in $X$ are scaled by a factor $t$ as one approaches $\partial {\mathbb B}$.

\section{Cone Differential Operators} 
Let $A$ be a differential operator of order $\mu$ on ${\mathbb B}$ with smooth coefficients up to the boundary. In local coordinates near $\partial {\mathbb B}$  it can be written in the form 
$$A = \sum_{j+|{\alpha}|\le \mu} c_{j,{\alpha}}(t,x) D_x^{\alpha} D_t^j$$
with smooth functions $c_{j,{\alpha}}$. Regrouping the coefficients and exchanging $D_t = -i\partial_t$ by $\partial_t$ we may rewrite $A$ in that neighborhood in the form 
\begin{eqnarray}
\label{eq.A}
A= t^{-\mu} \sum_{j=0} ^\mu a_j(t,x,D_x) (-t\partial_t)^j,
\end{eqnarray}
 where $t\mapsto a_j(t,x,D_x)$ is a smooth family of differential operators of order $\le\mu-j$ on $X$. 
We call operators that can be written close to $\partial {\mathbb B}$ in the form \eqref{eq.A} conically degenerate or cone differential operators. 
Note that the class of operators of the form \eqref{eq.A} is much larger than the class of smooth differential operators on ${\mathbb B}$.

\begin{example}
The Laplace-Beltrami operator $\Delta_g$ with respect to a metric $g$ that, in local coordinates is given by a matrix $(g_{jk}(y))$ and variables $y=(y_1, \ldots,y_m)$ is defined by the formula 
\begin{eqnarray}
\label{eq.LB}
{\Delta}_g = \frac1{\sqrt{\det g}} \sum_{j,k} \partial_{y_j} (g^{jk} \sqrt{\det g})\partial_{y_k}. 
\end{eqnarray}
where $(g^{jk})$ is the inverse matrix to $(g_{jk})$. 

It is instructive to compute the Laplace-Beltrami operator with respect to the metric $g$ in \eqref{eq.matrixg}.  Here, $\det g = t^{2n} \det h(t)$ and the inverse matrix to \eqref{eq.matrixg} is  
\begin{eqnarray}
\label{eq.ginv}
(g^{jk}) = 
\begin{pmatrix}
1&0 \ldots 0 \\
0& \\
\vdots &t^{-2}h^{-1}(t)\\
0&   
\end{pmatrix} .
\end{eqnarray}	
Denoting by   ${\Delta}_{h(t)} $  the Laplace-Beltrami operator with respect to the metric $h(t)$ on $X$ and letting $F(t) =\frac12 t\partial_t(\det h(t))/\det h(t)={\mathcal O}(t)$, we obtain
\begin{eqnarray}
\label{eq.LBg}
{\Delta}_g = t^{-2} \left( (-t\partial_t )^{2} -(n-1+F(t)) (-t\partial_t) + \Delta_{h(t)}\right).
\end{eqnarray}
Note that $F\equiv 0$ if the cone is straight, i.e. $h$ is independent of $t$.
\end{example}  
We see that $\Delta_g$ is of the form \eqref{eq.A} with $\mu=2$, $a_2(t,x,D_x) =1$, $a_1(t,x, D_x) =-(n-1+F(t))$ and $a_0(t,x,D_x)={\Delta}_{h(t)} $. 

\section{The Mellin transform} 
We denote by  $\Gamma_{\beta} =\{ z\in {\mathbb C}: {\mathop{\rm Re\, }} z={\beta}\}$ the vertical line of all complex numbers of real part ${\beta}$. 
\begin{definition}
Let $u\in C^\infty_c({\mathbb R}_{>0})$ (or, more generally, $u\in C^\infty_c({\mathbb R}_{>0},E)$ for some locally convex space $E$). The Mellin transform ${\mathcal M} u$ of $u$ is defined by 
$${\mathcal M} u(z) = \int_0^\infty t^z u(t) \frac{dt}t, \quad z\in {\mathbb C}.$$ 
\end{definition}
The following properties are easy consequences of the definition. 
\begin{lemma}
\label{4.2}
Let $u\in C^\infty_c({\mathbb R}_{>0})$.
\begin{enumerate}\renewcommand{\labelenumi}{(\alph{enumi})}%
\item ${\mathcal M} u$ is an entire function, i.e. holomorphic on ${\mathbb C}$.
\item ${\mathcal M}(t^{\gamma} u )(z) = ({\mathcal M} u)(z+{\gamma})$, ${\gamma}\in {\mathbb R}$. 
\item ${\mathcal M}(\ln t u )(z) = \partial_z{\mathcal M} u(z)$.
\item ${\mathcal M}(-t\partial_t u )(z) = z {\mathcal M} u(z)$. 
\item ${\mathcal M}$ extends to an isomorphism ${\mathcal M}: L^2({\mathbb R}_{>0} ) \to L^2(\Gamma_{1/2})$. 
\end{enumerate} 
\end{lemma} 

\begin{lemma}
\label{4,3}
 Let ${\omega}$ be a cut-off function on ${\mathbb R}_{>0}$, i.e. $0\le {\omega}\le 1$, ${\omega}$ is non-increasing, ${\omega}(t)\equiv 1$ near zero, ${\omega}(t)\equiv 0$ for large $t$. Then
$${\mathcal M}(t^{-p} \ln^kt{\omega}(t))(\cdot)$$
is meromorphic on ${\mathbb C}$ with a single pole of order $k+1$ in $z=p$. Outside any neighborhood of $p$, it is rapidly decreasing on each line $\Gamma_{\beta}$, uniformly for ${\beta}$ in compact intervals. 
\end{lemma} 

\subsection*{The weighted Mellin transform}
According to Lemma \ref{4.2}(b) we have 
\begin{eqnarray*}
({\mathcal M} u)_{|\Gamma_{1/2-{\gamma}}} (z) = {\mathcal M}(t^{-{\gamma}} u)(z+{\gamma}), \quad u\in C_c^\infty({\mathbb R}_{>0}).
\end{eqnarray*}

Since ${\mathcal M}: L^2({\mathbb R_{>0}})\stackrel\cong\to L^2(\Gamma_{1/2})$, this yields an extension of  ${\mathcal M}$ to an isomorphism
\begin{eqnarray*}
 {\mathcal M}_{\gamma}: t^{\gamma} L^2({\mathbb R_{>0}}) \to L^2(\Gamma_{1/2-{\gamma}}).
\end{eqnarray*}
The inverse is given by 
\begin{eqnarray*}
{\mathcal M}_{\gamma}^{-1} (t) =\frac1{2\pi i}  \int_{\Gamma_{1/2-{\gamma}}} t^{-z} h(z) \, dz. 
\end{eqnarray*}

Using Lemma \ref{4.2} and the inversion formula for the Melllin transform, we can write the operator $A$ with the representation \eqref{eq.A} in the form 
\begin{eqnarray}\label{eq.opm}
Au(t,x) = \frac1{2\pi i}\int_{\Gamma_{\frac{1}2-\gamma}}\!\!\! t^{-z} \sum_{j=0}^\mu a_j(t,x,D_x)z^j({\mathcal M}_{t\to z} u)(z,x)\, dz, \quad u\in C_c^\infty({\mathbb B}^\circ),
\end{eqnarray}
for $\gamma\in {\mathbb R}$. In fact, an application of Cauchy's integral formula shows that the result is the same for all ${\gamma}$, since  $u\in C^\infty_c({\mathbb B}^\circ)$, so that ${\mathcal M} u$ is holomorphic in $z$ and rapidly decreasing along vertical lines $\Gamma_{{\beta}}$, uniformly for ${\beta}$ in compact intervals. 

\section{Cone Sobolev Spaces} 
The form \eqref{eq.A} of a cone differential operator suggests that one should consider these operators not on standard Sobolev spaces but rather on a scale of spaces that 
\begin{itemize}
\item incorporates the weight $t^{-\mu}$ and, more generally, power weights near $t=0$ 
\item measures smoothness not with respect to derivatives $\partial_{x_j}$, $j=1,\ldots,n$, and $\partial_t$ but rather with respect to $\partial_{x_j}$ and $t\partial_t$.
\end{itemize}

\begin{definition}
For $s=k\in {\mathbb N}_0,{\gamma}\in {\mathbb R}$, $1<p<\infty$ let
\begin{eqnarray*}
{\mathcal H}^{k,{\gamma}}_p({\mathbb B}) &=&\left\{u\in H^k_{loc} ({\mathbb B}^\circ) :\right. \\
&&\left.t^{\frac{n+1}2-{\gamma}}\partial^{\alpha}_x(t\partial_t)^j u(t,x){\omega}(t)\in L^p\left({\mathbb B}, \frac{dxdt}t\right), \forall j+|{\alpha}|\le k\right\}.
\end{eqnarray*} 
Here, ${\omega}$ is an arbitrary cut-off function in $C^\infty_c(\overline{\mathbb R}_{>0})$.

This definition extends to $s\in {\mathbb R}$ by interpolation and duality. Note that $s$ measures smoothness and ${\gamma}$ flatness near $t=0$. 
\end{definition}

\begin{remark}
The weight and the measure in the definition of the spaces ${\mathcal H}^{s,{\gamma}}_p({\mathbb B})$ are motivated by the fact that in this way, ${\mathcal H}^{0,0}_2({\mathbb B})= L^2({\mathbb B},g)$, the $L^2$-space with respect to the  metric $g$. In fact, as we saw, $\sqrt{\det g} =t^n\sqrt{\det h(t)}$. Since $\sqrt{\det h(t)}$ is bounded and bounded away from zero we see that, for $u$ supported in $[0,1)\times X$, 
\begin{eqnarray*}
u\in L^2({\mathbb B},g) &\Leftrightarrow& \int _0^1 \int_X |u(t, x) |^2 t^n\sqrt{\det h(t)} \, dtdx<\infty\\
&\Leftrightarrow& t^{\frac{n+1}2} u(t, x)  \in L^2\Big(\frac{dtdx}t\Big) 
\Leftrightarrow u\in {\mathcal H}^{0,0}_2({\mathbb B}).  
\end{eqnarray*}
Note that this equality does not extend to $p\not=2$. We instead have $L^p({\mathbb B},g) = {\mathcal H}^{0,{\gamma}_p} ({\mathbb B})$, ${\gamma}_p=(n+1)(\frac12-\frac1p)$.
The advantage of the above definition is that, in this way, the spectrum of an operator $A$ is independent of $p$, see \cite{SS2001}. 
\end{remark} 

\begin{lemma} \begin{enumerate}\renewcommand{\labelenumi}{(\alph{enumi})}%
\item ${\mathcal H}^{s,{\gamma}}_p({\mathbb B}) ' = {\mathcal H}^{-s,-{\gamma}}_{p'}({\mathbb B})$ , where $1/p+1/p'=1$. 
\item ${\mathcal H}^{s,{\gamma}}_p({\mathbb B}) \hookrightarrow {\mathcal H}^{s',{\gamma}'}_p({\mathbb B})$ if and only if $s\ge s'$ and $ {\gamma}\ge{\gamma}'$.
\item ${\mathcal H}^{s,{\gamma}}_p({\mathbb B}) \stackrel{\text{\rm comp}}\hookrightarrow {\mathcal H}^{s',{\gamma}'}_p({\mathbb B})$ if and only if $s> s'$ and $ {\gamma}>{\gamma}'$.
\end{enumerate}
Hence one has to gain smoothness and weight in order to obtain compactness. 
\end{lemma}

For cone differential operators the following theorem is obvious. 
\begin{theorem}
Let $A$ be a conically degenerate operator of order $\mu$. Then 
$$ A:{\mathcal H}^{s+\mu,{\gamma}+\mu}_p({\mathbb B}) \longrightarrow {\mathcal H}^{s,{\gamma}}_p({\mathbb B})$$
is bounded. 
\end{theorem}

\section{Symbols} 
Recall that a differential operator on a smooth manifold $M$ has an invariantly defined principal symbol on $T^*M$. Namely, if $P=\sum_{|{\alpha}|\le \mu} c_{\alpha}(y) D^{\alpha}_y$, then 
$$
{\sigma}_\psi^\mu(P)(x,\xi) = \sum_{|{\alpha}|=\mu} c_{\alpha}(y)\xi^{\alpha}. 
$$
Choosing coordinates $(t,x,\tau,\xi)$ for $T^*([0,1)\times X)$ and recalling   that $D_t=-i\partial_t$, we obtain for a cone differential operator of the form \eqref{eq.A}
$${\sigma}_\psi^\mu(A)(t,x,\tau,\xi) = t^{-\mu} \sum_{j=0}^\mu {\sigma}_\psi^{\mu-j} (a_j) (t,x,\xi) (-it\tau)^j.$$
We notice the degeneracy in $t$ as $t\to0$ that is inherent to the singularity and therefore should not play a role when analyzing objects in the conical setting. We define the {\em rescaled} symbol as
$$t^\mu {\sigma}^\mu_\psi(A)(t, x,\tau/t,\xi) = \sum_{j=0}^\mu  {\sigma}_\psi^{\mu-j} (a_j) (t,x,\xi) (-i\tau)^j.$$

\begin{example}For the Laplace-Beltrami operator ${\Delta}_g$ Equation \eqref{eq.LBg} implies that 
$${\sigma}_\psi^2({\Delta}_g) (t, x,\tau,\xi) = t^{-2} ((-it\tau)^2 +{\sigma}_\psi^2({\Delta}_{h(t)})(x,\xi) ).$$
Its rescaled symbol therefore is 
$$t^2{\sigma}_\psi^\mu(\Delta_g)(t,x,\tau/t,\xi) = \tau^2 + {\sigma}^2_\psi({\Delta}_{h(t)})(x,\xi) = \tau^2 + \|\xi\|_{h^*(t)}^2,$$
where 
$ \|\xi\|_{h^*(t)}^2= \sum h^{jk}(t,x) \xi_j\xi_k$ is the squared norm of $\xi$ with respect to the metric $h^*(t) $ on the cotangent space induced from the metric $h(t)$. 
\end{example} 

For conically degenerate operators a second symbol is important. The so-called conormal symbol or principal Mellin symbol of a cone differential operator $A$, denoted by ${\sigma}_M^\mu(A)$  is a function on ${\mathbb C}$ with values in differential operators on $X$. For $A$ as in \eqref{eq.A}, we define
\begin{eqnarray*}
({\sigma}_M^\mu(A))(z)  = \sum_{j=0}^\mu a_j(0,x,D_x)z^j, \quad z\in {\mathbb C}.
\end{eqnarray*}
Hence we obtain bounded maps
$${\sigma}_M^\mu(A)(z) : H^{s+\mu} (X) \to H^s(X), \quad z\in {\mathbb C}.$$

\begin{definition} \label{6.2} 
A cone differential operator $A$ with a representation as in \eqref{eq.A} is called elliptic with respect to the weight line $\Gamma_{\beta}$, if 
\begin{enumerate}\renewcommand{\labelenumi}{(\roman{enumi})}
\item The principal pseudodifferential symbol ${\sigma}_\psi^\mu$ is invertible in the interior of ${\mathbb B}$ and, near $\partial {\mathbb B}$, the rescaled principal symbol is uniformly invertible up to $t=0$, i.e. its inverse is ${\mathcal O}(|\xi|+|\tau|)^{-\mu} ) $ up to $t=0$.
\item  ${\sigma}_M^\mu(A)(z) : H^{s+\mu} (X) \to H^s(X)$ is invertible for some $s\in {\mathbb R}$ and all $z\in \Gamma_{\beta}$. 
\end{enumerate} 
\end{definition} 
A standard result from the theory of pseudodifferential operators implies that ${\sigma}_M^\mu(A)(z) $ in (ii) is invertible for {\em all} $s\in {\mathbb R}$ if it is invertible for some $s$. 

\begin{example}
For the Laplace-Beltrami operator 
$$\Delta_g = t^{-2}\left((-t\partial_t)^2 - (n-1-F(t))(-t\partial_t) +{\Delta}_{h(t)}\right)$$ 
the principal Mellin symbol is {\rm(}recall that $F(t)={\mathcal O}(t)${\rm)}
$${\sigma}^2_M({\Delta}_g)(z) = z^2 -(n-1)z +{\Delta}_{h(0)}.$$
\end{example} 
For which $z$ is it invertible? This is the case if and only if $z^2-(n-1)z$ is not in the spectrum of $-{\Delta}_{h(0)}$. Denoting by $0={\lambda}_0>{\lambda}_1>\ldots $ the {\em different} eigenvalues of ${\Delta}_{h(0)}$ we see that ${\sigma}_M^2({\Delta}_g)(z)$ is {\em not} invertible, if, and only if 
$$z\in \Big\{ q_j^\pm = \frac{n-1}2\pm \sqrt{\Big(\frac{n-1}2\Big)^2-{\lambda}_j}, j\in {\mathbb N}_0 \Big\}.$$
In particular, ${\sigma}_M^2(\Delta_g)(\cdot)$ is invertible on every line $\Gamma_{{\beta}}$ with ${\beta}\not= q^\pm_j$, $j\in{\mathbb N}_0$.
As $X$ is a closed manifold, ${\lambda}_0=0$ and so ${\sigma}_M^2({\Delta}_g)(z)$ is not invertible in  $z=q_0^-=0$ and $z=q_0^+=n-1$ (both coincide for $n=1$).

\section{The Cone Calculus} 

Let  ${\omega}_1, {\omega}_2,{\omega}_3$ be cut-off functions  
with ${\omega}_1{\omega}_2 = {\omega}_1$, ${\omega}_1{\omega}_3={\omega}_3$ and ${\gamma}\in {\mathbb R}$. 

An operator of order $\mu \in {\mathbb R}$ in the cone calculus with respect to the line $\Gamma_{\frac{n+1}2-{\gamma}}$  is an operator $A: C^\infty_c({\mathbb B}^\circ) \to C^\infty({\mathbb B}^\circ)$ that can be written in the form (cf.\eqref{eq.opm})
\begin{eqnarray}
\label{eq.cA}
A= t^{-\mu} {\omega}_1 {\text{\rm op}}_M^{\gamma}(h+h_0) {\omega}_2 + (1-{\omega}_1)P(1-{\omega}_3) + G,
\end{eqnarray}
where
\begin{enumerate}\renewcommand{\labelenumi}{(\roman{enumi})} 
\item $h=h(t,z)$ is  smooth in $t\ge0$ and holomorphic in $z\in {\mathbb C}$, taking values in pseudodifferential operators of order $\mu$ on $X$ such that $h(t, z)_{|z\in \Gamma_{{\beta}}}\in \Psi^\mu(X;\Gamma_{{\beta}}) $ 
is a parameter-dependent (w.r.t. $\Gamma_{\beta}$) pseudodifferential operator of order $\mu$, uniformly for ${\beta}$ in compact intervals.
\item $h_0=h_0(z)$ is meromorphic on ${\mathbb C}$ with values in pseudodifferential operators of order $\mu$ on $X$. The coefficients in the  Laurent expansions in the poles are pseudodifferential operators of finite rank, so in particular smoothing. No pole lies on $\Gamma_{({n+1)}/2-{\gamma}}$ and 
$h_{0|\Gamma_{{\beta}}}\in \Psi^\mu(X;\Gamma_{{\beta}})$ outside the poles, unformly for ${\beta}$ in compact intervals.  
\item $P$ is a pseudodifferential operator of order $\mu$ on ${\mathbb B}^\circ$.  
\item $G:{\mathcal H}^{s,{\gamma}} ({\mathbb B})\to {\mathcal H}^{\infty,\gamma-\mu+{\varepsilon}}({\mathbb B})$ and $G^*: {\mathcal H}^{s-+\nu,-{\gamma}+\mu}({\mathbb B}) \to {\mathcal H}^{\infty,-{\gamma}+{\varepsilon}}({\mathbb B})$  for all $s\in {\mathbb R}$ and  some ${\varepsilon}>0$.  
\end{enumerate}   

The operator ${\text{\rm op}}_M^{\gamma}(h+h_0):C^\infty_c({\mathbb B}^\circ) \to C^\infty(\mathbb B^\circ)$  is then defined by 
$${\text{\rm op}}_M^{\gamma}(h+h_0) u(t,\cdot) = \frac{1}{2\pi i}\int_{\Gamma_{\frac{n+1}2-{\gamma}}} t^{-z} (h(t,z )+h_0(z)) {\mathcal M}_{t\to z}  u(z,\cdot)dz.
$$
More details can be found in \cite{SS1994} and \cite{SS1995} or \cite[Chapter 2]{S1998}. 

The operators $G$ in (iv) form the residual class. They are infinitely smoothing and gain a little weight, so they  are compact elements of  ${\mathcal L}({\mathcal H}^{s,{\gamma}}({\mathbb B}), {\mathcal H}^{s-\mu,\gamma-\mu}({\mathbb B}))$.  

\begin{theorem}%
Let $A$ be as above. For every $s\in {\mathbb R}$,  $A$ in \eqref{eq.cA} then extends to a bounded linear operator 
$$ A: {\mathcal H}^{s,{\gamma}}({\mathbb B}) \to {\mathcal H}^{s-\mu,{\gamma}-\mu}({\mathbb B}).$$
\end{theorem} 

\begin{theorem} {\rm (a)}
The operators in the cone calculus form a graded $*$-algebra: 
The composition of a cone operator of order $\mu_1$ with respect to the line $\Gamma_{\frac{n+1}2-{\gamma}-\mu_2}$  with one of order $\mu_2$ w.r.t. $\Gamma_{\frac{n+1}2-{\gamma}}$ is an element of order $ \mu_1+\mu_2$ 
w.r.t. $\Gamma_{\frac{n+1}2-{\gamma}}$.

 {\rm (b)} The adjoint with respect to an extension of the $L^2$ inner product of an operator of order $\mu$  w.r.t. the line $\Gamma_{\frac{n+1}2-{\gamma}}$ is an operator of order $\mu$ w.r.t. $\Gamma_{\frac{n+1}2-{\gamma}-\mu}$.	\end{theorem} 

\subsection*{Ellipticity} 
The notions of ellipticity extend to the general cone calculus. In fact, to the operator ${\text{\rm op}}_M^{\gamma}(h)$ we may associate locally a conically degenerate pseudodifferential symbol of the form $p(t,x,t \tau, \xi) $, where $p(t, x,\tau, \xi)$ is a standard symbol of order $\mu$, smooth up to $t=0$, modulo symbols of order $-\infty$ in the interior. 

For the operator \eqref{eq.cA} we therefore have a well-defined principal symbol in the interior with conical degeneration as $t\to 0$. We can define the rescaled symbol as before by omitting the factor $t^{-\mu}$ and replacing $t\tau$ by $\tau$. 
As in Definition \ref{6.2}(a) ellipticity  requires the invertibility of the interior pseudodifferential symbol together with the uniform invertibility of the rescaled symbol.

Following \ref{6.2}(b), we define the principal Mellin symbol ${\sigma}^\mu_M(A) (z) = h(0,z)+h_0(z)$. This is a pseudodifferential operator of order $\mu$ on $X$. As in Definition \ref{6.2}, ellipticity with respect to the line $\Gamma_{\frac{n+1}2-{\gamma}}$ requires additionally the invertibility of ${\sigma}^\mu_M(A)(z): H^s(X) \to H^{s-\mu}(X)$ for some (and then all) $s\in {\mathbb R}$ and $z$ on this line. 

A central theorem then is: 

\begin{theorem} 
Assume that the cone operator $A$ is elliptic with respect to the line $\Gamma_{\frac{n+1}2-{\gamma}}$.  
Then there exists a parametrix $B$ of order $-\mu$ in the cone calculus such that 
$$AB-I_{{\mathcal H}_p^{s-\mu,\gamma-\mu}({\mathbb B}) }=G_R \text{ and } BA-I_{{\mathcal H}_p^{s,{\gamma}}({\mathbb B})}=G_L$$ 
are bounded operators $G_R: {\mathcal H}_p^{s-\mu,{\gamma}-\mu} ({\mathbb B})\to {\mathcal H}_p^{\infty,\gamma-\mu+{\varepsilon}}({\mathbb B})$ and 
$G_L: {\mathcal H}_p^{s,{\gamma}} ({\mathbb B}) \to {\mathcal H}_p^{\infty,\gamma+{\varepsilon}}({\mathbb B})$. 
In particular, $A$ then is a Fredholm operator. 
\end{theorem} 

\subsection*{Back to the beginning} One can now see the reason for the appearance of the singular solutions: In analogy to the case of pseudodifferential parametrices, $B$ is constructed starting from a rough version $B_0$ defined by inverting the principal pseudodifferential symbol and the principal Mellin symbol and subsequent refinements.  Here, 
$$B_0 =   \tilde {\omega}_1  {\text{\rm op}}_M^{\gamma}(g(z))\tilde  {\omega}_2t^\mu  + (1-\tilde {\omega}_1)Q(1-\tilde {\omega}_3)$$
with $g(z) = (h(0,z)+h_0(z))^{-1} $, a parametrix $Q$ to $P$ and cut-off functions $\tilde {\omega}_1, \tilde {\omega}_2,\tilde {\omega}_3$. 
According to Lemma \ref{4.2}(b)  we have for $u\in C^\infty_c({\mathbb B}^\circ)$
\begin{eqnarray*}
{\text{\rm op}}_M^{\gamma} (g)(t^\mu u)(t,x) &=& \frac1{2\pi i} \int_{\Gamma_{(n+1)/2-{\gamma}}} t^{-z} g(z) {\mathcal M} (t^{\gamma} u)(z,x) \, dz\\
&=&  \frac1{2\pi i} \int_{\Gamma_{(n+1)/2-{\gamma}}} t^{-z} g(z) {\mathcal M} (u)(z+\mu,x) \, dz\\
&=&  \frac1{2\pi i} \int_{\Gamma_{(n+1)/2-{\gamma}+\mu}} t^{-w+\mu} g(w-\mu) {\mathcal M} (u)(w,x) \, dw\\
&=& t^\mu {\text{\rm op}}_M^{\gamma-\mu} (g(\cdot - \mu)) u(t,x)\\
&=& t^\mu {\text{\rm op}}_M^{\gamma} (g(\cdot-\mu)) u(t, x) + \frac1{2\pi i} \int_{\mathcal C} t^{-w} g(w) {\mathcal M} (u)(w,x) \, dw,
\end{eqnarray*}
where $\mathcal C$ is the `contour' in ${\mathbb C}$ given by the vertical lines ${\Gamma_{(n+1)/2-{\gamma}}}$ and ${\Gamma_{(n+1)/2-{\gamma}+\mu}}$. The integral makes sense, since $g$ decreases rapidly at infinity, and the integration furnishes the residues of the poles of $w\mapsto t^{-w} g(w-\mu){\mathcal M} u(w,x)$ between the lines. We know that ${\mathcal M} u(w,\cdot)$ is holomorphic in $w$ and the Laurent coefficients of $g$ in a pole are finite rank (hence smoothing) pseudodifferential operators. So a pole of $g$ in the point $p$ between the two lines of order $k\in {\mathbb N}$ furnishes contributions of the form $t^{-p} \ln^{j}(t) s_j(x)$, where the $s_j$ are smooth functions on $X$, $j\le k-1$. Subsequent refinements lead additionally to  terms  $t^{-p+l} \ln^{j}(t) \tilde s_j(x)$, $\tilde s_j\in C^\infty(X)$, $l\in {\mathbb N}$, as long as $p-l$ also lies between the lines.

\section{Historical Remarks}Kondrat'ev's paper \cite{K1967} from 1967 is generally considered as the starting point of the analysis on manifolds with conical singularities. His approach was further developed e.g. by V.G.  Maz'ja and   B.A.  Plamenevski\u \i\ \cite{MP1973}. Independently, Jeff Cheeger in \cite{C1979} studied  the spectral geometry of spaces with cone-like singularities.  However, it was only in the 1980's that pseudodifferential calculi were systematically developed. 
Early works are \cite{MM1983} by Richard Melrose and Gerardo Mendoza, \cite{S1986} by B.-W. Schulze and the monographs \cite{P1990} by  Plamenevski\u \i, \cite{S1991} by Schulze and \cite{M1993} by Melrose. Several other singular calculi were developed in the sequel, in particular the edge calculus by Schulze, see e.g. \cite{RS1987}, and Mazzeo \cite{M1990}, the zero calculus and the $\phi$-calculus by Mazzeo and Melrose \cite{MM1987}, \cite{MM1998}, and the cone calculus for boundary value problems by  Schulze and the author \cite{SS1994}, \cite{SS1995}.   
Index theory on singular spaces was studied by many authors, two examples are Br\"uning and Seeley \cite{BS1988} and Lesch \cite{L1997}. 
In collaboration with Sandro Coriasco, Nikolaos Roidos and J\"org Seiler, the author has developed an approach to the analysis of nonlinear evolution equations  on manifolds with conical singularities via maximal regularity techniques, starting with \cite{CSS2003} and relying  on results from \cite{SS2005}. The Cahn-Hilliard equation and the Allen-Cahn equation were treated in \cite{RS2013} and \cite{RS2014}, the porous medium equation in \cite{RS2016}, \cite{RS2018} and \cite{RSS2021}. 
}

\end{document}